\documentclass[12pt]{article}
\usepackage{latexsym,amsfonts,amssymb,amsmath}
\usepackage{color}
\usepackage{hyperref}
\usepackage{cleveref}
\usepackage{fancyhdr}
\usepackage{setspace}
\usepackage{colordvi,multicol}

\newtheorem{thm}{Theorem}[section]

\newtheorem{pro}[thm]{Proposition}

\newtheorem{rem}[thm]{Remark}

\date{}

\begin{document}

\title{\bf The infimum values of the probability functions for some infinitely divisible distributions  motivated by Chv\'{a}tal's theorem}

 \author{ Ze-Chun Hu$^{1}$, Peng Lu$^1$, Qian-Qian Zhou$^{2,}$\footnote{Corresponding author: qianqzhou@yeah.net}, Xing-Wang Zhou$^{1}$\\ \\
   {\small $^1$ College of Mathematics, Sichuan  University,
 Chengdu 610065, China}\\
 {\small $^2$ School of Science, Nanjing University of Posts and Telecommunications, Nanjing  210023, China}}
\maketitle

\makeatletter

\begin{abstract}

  Let $B(n,p)$ denote a binomial random variable with parameters $n$ and $p$.  Chv\'{a}tal's theorem  says that for any fixed $n\geq 2$, as $m$ ranges over $\{0,\ldots,n\}$, the probability $q_m:=P(B(n,m/n)\leq m)$ is the smallest when  $m$ is closest to $\frac{2n}{3}$.  Motivated by this theorem, in this paper  we consider the infimum value of the probability $P(X\leq \kappa E[X])$, where $\kappa$ is a positive real number, and  $X$ is  a random variable whose distribution belongs to some infinitely divisible distributions including  the inverse Gaussian, log-normal, Gumbel and logistic distributions.
\end{abstract}



\section{Introduction}

Let $B(n,p)$ denote a binomial random variable with parameters $n$ and $p$. Janson \cite{Ja21} introduced the following conjecture suggested by Va\v{s}k Chv\'{a}tal.

{\bf Conjecture 1} (Chv\'{a}tal). For any fixed $n\geq 2$, as $m$ ranges over $\{0,\ldots,n\}$, the probability $q_m:=P(B(n,m/n)\leq m)$ is the smallest when
 $m$ is closest to $\frac{2n}{3}$.

Chv\'{a}tal's conjecture has  applications in machine learning (see Doerr \cite{Do18}, Greenberg and Mohri \cite{GM14} and the references therein). Janson \cite{Ja21} showed that  Conjecture 1 holds for large $n$. Barabesi et al. \cite{BPR23} and Sun \cite{Su21} proved that  Conjecture 1 is true  for general $n\geq 2$. Hereafter, we call Conjecture 1 by Chv\'{a}tal's theorem.

Motivated by Chv\'{a}tal's theorem,  Li et al. \cite{LXH23}  considered the infimum value problem on the probability that a random variable is not more than its expectation, when its distribution is the Poisson distribution, the geometric distribution or the Pascal distribution, where as to the Pascal distribution, only some partial results have been obtained. Sun et al. \cite{SHS23} investigated  the corresponding infimum value problem for the Gamma distribution among other things.  Li et al. \cite{LHZ23} studied  the infimum value problem for the Weibull distribution and the Pareto distribution. Guo et al. \cite{GTH23} considered the infimum value of the probability $P(X\leq  E[X])$, where   $X$ is  a negative binomial random variable, and gave an affirmative answer to the conjecture on the Pascal distribution posed in \cite{LXH23}.

In this paper,  we consider the infimum value of the probability $P(X\leq \kappa E[X])$, where $\kappa$ is a positive real number, and $X$ is  a random variable whose distribution belongs to some infinitely divisible distributions including  the inverse Gaussian, log-normal, Gumbel and logistic distributions in Section 2 and Section 3.  Before presenting the main results, we give a remark.

{\bf Remark 1.1}\ Let $X_{\alpha}$ be a random variable with finite expectation $E[X_{\alpha}]$,  where $\alpha$ stands for some parameter in the distribution of $X_{\alpha}$, and it is a real number or a vector in $\mathbf{R}^n$ for some positive integer $n\geq 2$. We have the following two motivations to study $\inf_{\alpha}P(X_{\alpha}\leq \kappa E[X_{\alpha}])$:

\begin{itemize}

\item From $\inf_{\alpha}P(X_{\alpha}\leq   \kappa E[X_{\alpha}])$, we can get $\sup_{\alpha}P(X_{\alpha}>\kappa E[X_{\alpha}])$. Obviously, if we wish the probability $P(X_{\alpha}>\kappa E[X_{\alpha}])$ is as large as possible, we should consider $\sup_{\alpha}P(X_{\alpha}>\kappa E[X_{\alpha}])$ or equivalently $\inf_{\alpha}P(X_{\alpha}\leq \kappa E[X_{\alpha}])$. Based on this observation, we wish that our work on this topic may find some applications in machine learning, statistics, finance and economics etc.

\item Assume that $X_{\alpha}$ is nonnegative and denote by $\mu_{\alpha}$ the distribution of $X_{\alpha}$. If $\inf_{\alpha}P(X_{\alpha}\leq \kappa E[X_{\alpha}])=\beta>0$, then for any $\alpha$, we have
    \begin{align*}
\mu_{\alpha}([0,\kappa E[X_{\alpha}])\geq \beta,
    \end{align*}
which tells us that the family of the distributions $\{\mu_{\alpha}\}$ possesses a kind of  measure concentration phenomenon.

\end{itemize}

\section{Inverse Gaussian distribution}\setcounter{equation}{0}

Let $X_{\mu,\lambda}$  be an inverse Gaussian random variable with parameters $\mu$  and $\lambda$ $(\mu>0,\lambda>0).$ By \cite[Chapter 27]{K}, we  know that
the density function of $X_{\mu,\lambda}$  is given by
$$
f_{\mu , \lambda}(x)=\sqrt{\frac {\lambda}{2\pi x^3}}\exp\left(-\frac{\lambda(x-\mu)^2}{2{\mu}^2x}\right),\  \  x>0,
$$
and $E[X_{\mu,\lambda}]=\mu$.  Then, for any given real number $\kappa>0$, by \cite[Chapter 27]{K},  we know that
\begin{eqnarray*}
P(X_{\mu,\lambda}\leq \kappa E[X_{\mu,\lambda}])&=& \int_0^{\kappa \mu} \sqrt{\frac {\lambda}{2\pi x^3}}\exp\left(-\frac{\lambda(x-\mu)^2}{2{\mu}^2x}\right) dx\\
&=&\Phi\left(\sqrt{\frac{\lambda}{\kappa \mu }}(\kappa-1)\right)+\exp\left(\frac{2\lambda}{\mu}\right)\Phi\left(-\sqrt{\frac{\lambda}{\kappa \mu}}(\kappa+1)\right).
\end{eqnarray*}
Hereafter, $\Phi(\cdot)$ stands for the distribution function of the standard normal distribution.

Denote $x:=\sqrt{\frac{\lambda}{\mu}}$.  Then
$$
P(X_{\mu,\lambda}\leq \kappa E[X_{\mu,\lambda}])=\Phi\left(\sqrt{\frac{1}{\kappa}}(\kappa-1)x\right)+e^{2x^2}
\Phi\left(-\sqrt{\frac{1}{\kappa}}(\kappa+1)x\right).
$$
Define a function
\begin{gather}\label{d1}
g_{\kappa}(x):=\Phi\left(\sqrt{\frac{1}{\kappa}}(\kappa-1)x\right)+e^{2x^2}
\Phi\left(-\sqrt{\frac{1}{\kappa}}(\kappa+1)x\right),\  x>0.
\end{gather}

The main result of this section is
\begin{thm}\label{pro-4.1}
$(i)$ If $\kappa \leq 1$, then
\begin{eqnarray*}
\inf\limits_{x \in \left ( 0,\infty  \right ) }g_{\kappa }\left ( x \right )=\lim_{x\to \infty}g_{\kappa}(x)=
\left\{
\begin{array}{cl}
0,& \mbox{if}\ \kappa < 1,\\
\frac{1}{2},& \mbox{if}\ \kappa = 1.
\end{array}
\right.
\end{eqnarray*}
$(ii)$ If $\kappa>1$, then
\begin{eqnarray*}
\min_{x\in \left ( 0,\infty  \right ) }g_{\kappa }\left ( x \right )=g_{\kappa }\left ( x_{0}(\kappa) \right ),
\end{eqnarray*}
where $x_{0}(\kappa)$ is the unique zero point of the function
\begin{eqnarray}\label{2.2}
h_{\kappa}(x):=2\int^{\infty}_{\sqrt{\frac{1}{\kappa}}(\kappa+1)x}
\exp\left(-\frac{1}{2}t^2\right)dt-\sqrt{\frac{1}{\kappa}}\frac{\exp(-\frac{{(\kappa+1)}^2}
{2\kappa}x^2)}{x},\  x\in (0, \infty).
\end{eqnarray}
\end{thm}

\noindent {\bf Proof. } By taking the derivative of the function $g_{\kappa}(x)$  defined in \eqref{d1} and in virtue of (\ref{2.2}), we have
\begin{eqnarray}\label{2.3}
g_{\kappa }^{'}(x)&=&\sqrt{\frac{1}{\kappa}}(\kappa-1)\Phi'\left(\sqrt{\frac{1}{\kappa}}(\kappa-1)x\right)
+4xe^{2x^2}
\Phi\left(-\sqrt{\frac{1}{\kappa}}(\kappa+1)x\right)\nonumber\\
&&-\sqrt{\frac{1}{\kappa}}(\kappa+1)e^{2x^2}
\Phi'\left(-\sqrt{\frac{1}{\kappa}}(\kappa+1)x\right)\nonumber\\
&=&\frac{4xe^{2x^2}}{\sqrt{2\pi}}\int_{\sqrt{\frac{1}{\kappa}}(\kappa+1)x}^{\infty}
\exp\left(-\frac{1}{2}t^2\right)dt-2\sqrt{\frac{1}{2\pi\kappa}}
\exp\left(-\frac{{(\kappa-1)}^2}{2\kappa}x^2\right)\nonumber\\
&=&\frac{2xe^{2x^2}}{\sqrt{2\pi}}\left[
2\int_{\sqrt{\frac{1}{\kappa}}(\kappa+1)x}^{\infty}\exp\left(-\frac{1}{2}t^2\right)dt-
\sqrt{\frac{1}{\kappa}}\frac{\exp\left(-\frac{{(\kappa+1)}^2}{2\kappa}x^2\right)}{x}\right]\nonumber\\
&=&\frac{2xe^{2x^2}}{\sqrt{2\pi}}h_{\kappa}(x).
\end{eqnarray}
By taking the derivative of the function $h_{\kappa}(x),$
 we get that
\begin{eqnarray}\label{d2}
h_{\kappa }^{'}(x)&=&\sqrt{\frac{1}{\kappa}}\exp\left(-\frac{(\kappa+1)^2}{2\kappa}x^2\right)
\left(\frac{1}{\kappa}-\kappa+\frac{1}{x^2}\right)\nonumber\\
&:=&\sqrt{\frac{1}{\kappa}}\exp\left(-\frac{(\kappa+1)^2}{2\kappa}x^2\right)
\varphi_{\kappa}(x),
\end{eqnarray}
where $\varphi_{\kappa}(x)=\frac{1}{\kappa}-\kappa+\frac{1}{x^2}.$

$(i)$ If $\kappa \leq 1$ , then  $\varphi_{\kappa}(x)\geq \frac{1}{x^2}>0$ for any    $ x \in \left (0 , \infty \right)$. Thus, we have $h_{\kappa }^{'}(x)>0,\ \  \forall  x\in (0, \infty),$  which implies  that the function $h_{\kappa}(x)$ is strictly increasing on interval $(0, \infty)$. We have
\begin{eqnarray*}
\lim_{x \rightarrow \infty} h_{\kappa }\left ( x \right )
& =& \lim_{x \rightarrow \infty}\left(2\int_{\sqrt{\frac{1}{\kappa}}(\kappa+1)x}^{\infty}
\exp\left(-\frac{1}{2}t^2\right)dt-
\sqrt{\frac{1}{\kappa}}\frac{\exp(-\frac{{(\kappa+1)}^2}{2\kappa}x^2)}{x}\right)=0.
\end{eqnarray*}
Thus for any $x\in (0,\infty)$,
$h_{\kappa }(x)< 0$. Then by (\ref{2.3}), we get
$$
g_{\kappa }^{'}\left ( x \right )
 < 0,\ \forall x\in (0,\infty).
 $$
Therefore, $g_{\kappa }\left ( x \right )$ is a strictly decreasing function on interval $ \left (0,\infty \right)$ and thus
\begin{eqnarray*}
&&\inf_{x\in (0, \infty)}g_{\kappa}(x)=\lim_{x\rightarrow \infty}g_{\kappa}(x)\\
&&=\lim_{x\rightarrow \infty}\Phi\left(\sqrt{\frac{1}{\kappa}}(\kappa-1)x\right)+\lim_{x\rightarrow \infty}e^{2x^2}\frac{1}{\sqrt{2\pi}}\int_{\sqrt{\frac{1}{\kappa}}(\kappa+1)x}^{\infty}
\exp\left(-\frac{t^2}{2}\right)dt\\
&&=\lim_{x\rightarrow \infty}\Phi\left(\sqrt{\frac{1}{\kappa}}(\kappa-1)x\right)+\lim_{x\rightarrow \infty}\frac{1}{\sqrt{2\pi}}e^{2x^2}\int_{0}^{\infty}
\exp\left(-\frac{{\left(u+\sqrt{\frac{1}{\kappa}}(\kappa+1)x\right)}^2}{2}\right)du\\
&&=\lim_{x\rightarrow \infty}\Phi\left(\sqrt{\frac{1}{\kappa}}(\kappa-1)x\right)+\lim_{x\rightarrow \infty}\frac{1}{\sqrt{2\pi}}\int_{0}^{\infty}\exp\left(-\frac{1}{2}\left(u^2
+\frac{2(\kappa+1)}{\sqrt {\kappa}}xu+\frac{({\kappa-1})^2}{\kappa}x^2\right)\right)du\\
&&=\lim_{x\rightarrow \infty}\Phi\left(\sqrt{\frac{1}{\kappa}}(\kappa-1)x\right)\\
&&=\left \{
\begin{array}{ll}
0,& \mbox{if}\ \kappa < 1,\\
\frac {1}{2},& \mbox{if}\ \kappa = 1.
\end{array}
\right.
\end{eqnarray*}

$(ii)$ If $\kappa > 1$, then by $h^{'}_{\kappa}(x)=0$, i.e., $\varphi_{\kappa}(x)=\frac{1}{\kappa}-\kappa+\frac{1}{x^2}=0$, we get $x=\sqrt{\frac{\kappa}{\kappa^2-1}}.$  Then, by \eqref{d2}, we  know  that
$h^{'}_{\kappa}(x)>0$  for any $x\in (0, \sqrt{\frac{\kappa}{\kappa^2-1}})$  and $h^{'}_{\kappa}(x)<0$ for any  $x\in (\sqrt{\frac{\kappa}{\kappa^2-1}}, \infty)$. It follows that   the function $h_{\kappa}(x)$  is strictly increasing on   $(0, \sqrt{\frac{\kappa}{\kappa^2-1}})$ and strictly decreasing on $(\sqrt{\frac{\kappa}{\kappa^2-1}}, \infty)$.  Furthermore, by  the fact that  $\lim_{x \rightarrow 0+}h_{\kappa}(x)=-\infty <0$ and $\lim_{x \rightarrow \infty}h_{\kappa}(x) = 0$, we know that  the continuous function $h_{\kappa}(x)$  has a unique zero point $x_{0}(\kappa)\in (0, \sqrt{\frac{\kappa}{\kappa^2-1}})$.
It follows that $h_{\kappa}(x)<0$ for any  $x\in (0, x_0(\kappa))$ and $h_{\kappa}(x)>0$  for any  $x\in (x_0(\kappa), \infty)$.

Hence, by \eqref{2.3}, we get that $g^{'}_{\kappa}(x)<0$  on $(0, x_0(\kappa))$  and $g^{'}_{\kappa}(x)>0$  on  $(x_0(\kappa), \infty)$. It follows that
$$\min_{x\in (0, \infty)}g_{\kappa}(x)=g_{\kappa}(x_0(\kappa)).$$ The proof is complete. \hfill\fbox

\begin {rem}
(i) For any $\kappa>1, \mu>0,\lambda>0$, it is easy to know that
$$
P(X_{\mu,\lambda}\le \kappa E[X_{\mu,\lambda}])>P(X_{\mu,\lambda}\le E[X_{\mu,\lambda}])>\frac{1}{2},
$$
where the second inequality follows from the proof of Theorem \ref{pro-4.1}. Then
for any $\kappa>1$, we have
\begin{eqnarray*}
\min_{x\in \left ( 0,\infty  \right ) }g_{\kappa }\left ( x \right )=g_{\kappa }\left ( x_{0}(\kappa) \right )>\frac{1}{2}.
\end{eqnarray*}

(ii) The above analysis shows that there is an interesting phase transition phenomenon
in the infimum value problem for the inverse Gaussian distribution. The critical point is $\kappa=1$.
\end {rem}

\section{Log-normal, Gumbel and Logistic distributions}

\subsection{Log-normal distribution}\setcounter{equation}{0}
Let $X_{\mu,\sigma}$  be a log-normal random variable with parameters $\mu$ and $\sigma$ $(\mu \in \mathbb{R}, \sigma>0).$ By \cite[Chapter 22]{K}, we know that the density function of $X_{\mu,\sigma}$ is
$$f_{\mu ,\sigma}(x)=\frac {1}{\sqrt {2\pi }\sigma x}\exp\left(-\frac {(\ln x-\mu)^2}{2{\sigma}^2}\right),\ \   x>0,$$
and   the expectation  is $E[X_{\mu,\sigma}]=\exp(\mu +\frac{{\sigma}^2}{2})$.
 Then,  for any given real number  $\kappa>0,$ by \cite[(22.1.2)]{K},  we have
\begin{eqnarray*}
P(X_{\mu,\sigma}\leq \kappa E[X_{\mu,\sigma}])=\Phi\left(\frac{\ln\left(\kappa e^{\mu+\frac{\sigma^2}{2}}\right)-\mu}{\sigma}\right)
=\Phi\left(\frac{\ln\kappa}{\sigma}+\frac{\sigma}{2}\right),
\end{eqnarray*}
which shows  that $P(X_{\mu,\sigma}\leq \kappa E[X_{\mu,\sigma}])$ is independent of $\mu.$

Define a function
$$
g_\kappa(\sigma):=\Phi\left(\frac{\ln\kappa}{\sigma}+\frac{\sigma}{2}\right),\   \sigma>0.
$$

The main result of this section is

\begin{pro}\label{pro-3.1}
$(i)$ If $\kappa \leq 1$, then
\begin{eqnarray}\label{Pro-3.1-a}
\inf\limits_{\sigma \in (0,\infty)}g_{\kappa }( \sigma)=\lim_{\sigma\to 0+}g_{\kappa}(\sigma)=
\left\{
\begin{array}{cl}
0,& \mbox{if}\ \kappa < 1,\\
\frac{1}{2},& \mbox{if}\ \kappa = 1;
\end{array}
\right.
\end{eqnarray}
$(ii)$ If $\kappa>1$, then
\begin{eqnarray}\label{Pro-3.1-b}
\min_{\sigma \in \left ( 0,\infty  \right ) }g_{\kappa }\left ( \sigma  \right )=g_{\kappa }\left ( \sqrt{2 \ln \kappa}  \right )=\Phi(\sqrt{2\ln \kappa })>\frac{1}{2}.
\end{eqnarray}
\end{pro}
{\bf Proof. } Define a function
\begin{gather}\label{c1}
h_{\kappa}(\sigma):=\frac{\ln \kappa}{\sigma}+\frac{\sigma}{2},\  \sigma>0.
\end{gather}
Then $g_{\kappa}(\sigma)=\Phi(h_{\kappa}(\sigma)).$ Therefore, in order to prove this proposition, it is enough to investigate  the value of $h_{\kappa}(\sigma)$.

 If $\kappa\le 1,$ then $\ln \kappa\le 0.$  Thus, by \eqref{c1},  we have $h^{'}_{\kappa}(\sigma)=-\frac{\ln \kappa}{\sigma^2}+\frac{1}{2}>0$, which implies that $h_{\kappa}(\sigma)$  is a strictly increasing function of $\sigma$.  Thus,
 \begin{eqnarray*}
\inf\limits_{\sigma \in \left ( 0,\infty  \right ) }h_{\kappa }\left ( \sigma  \right )=\lim_{\sigma\to 0+}h_{\kappa}(\alpha)=
\left\{
\begin{array}{cl}
-\infty,& \mbox{if}\ \kappa < 1,\\
0,& \mbox{if} \ \kappa = 1,
\end{array}
\right.
\end{eqnarray*}
It follows that (\ref{Pro-3.1-a}) holds.

 If $\kappa>1,$ then $\ln \kappa>0.$ If $h^{'}_{\kappa}(\sigma)=-\frac{\ln \kappa}{\sigma^2}+\frac{1}{2}=0$, then  $\sigma=\sqrt{ 2\ln \kappa}.$ It is easy to check that the function $h_{\kappa}(\sigma)$  is strictly decreasing on  $(0, \sqrt{ 2\ln \kappa})$ and strictly increasing on  $(\sqrt{ 2\ln \kappa}, \infty)$.  Thus,
 $$\min_{\sigma\in (0, \infty)}h_{\kappa}(\sigma)=h_{\kappa}(\sqrt{ 2\ln \kappa})=\sqrt{ 2\ln \kappa}.$$
 It follows that (\ref{Pro-3.1-b}) holds.   The proof is complete.  \hfill\fbox

\subsection{Gumbel distribution}\setcounter{equation}{0}

Let $X_{\mu,\beta}$ be a Gumbel random variable with parameters $\mu$ and $\beta$ $(\mu\in \mathbb{R},  \beta>0)$. By \cite{S73}, the density function of $X_{\mu,\beta}$ is given by
$$
f_{\mu,\beta}(x)=\frac{1}{\beta}e^{-z-e^{-z}},\  z=\frac{x-\mu}{\beta},\  x\in \mathbb{R}
$$
and its  expectation   is $E[X_{\mu,\beta}]=\mu +\beta \gamma,$  where $\gamma$  is  the Euler's constant. Then, for any given real number $\kappa>0,$  we have
$$
P(X_{\mu,\beta}\leq \kappa E[X_{\mu,\beta}])=e^{-e^{-\frac{\kappa{(\mu+\beta \gamma)}-\mu}{\beta}}}=e^{-e^{-\frac{(\kappa-1)\mu}{\beta}-\kappa \gamma}}.
$$

Let $x:=\frac{\mu}{\beta}$. Then
$$
P(X_{\mu,\beta}\leq \kappa E[X_{\mu,\beta}])=e^{-e^{-[(\kappa-1)x+\kappa \gamma]}}.
$$
Define a function
\begin{eqnarray}\label{4.1-a}
g_{\kappa}(x):=e^{-e^{-[(\kappa-1)x+\kappa \gamma]}},\  \   x\in \mathbb{R}.
\end{eqnarray}

The main result of this section is

\begin{pro}\label{pro-5.1}
$(i)$ If $\kappa <1$, then
\begin{eqnarray*}
\inf\limits_{x \in \mathbb{R}  }g_{\kappa }\left (x  \right )=\lim_{x\to \infty}g_{\kappa}(x)=0;
\end{eqnarray*}
$(ii)$ If $\kappa> 1$, then
$$
\inf_{x \in \mathbb{R} }g_{\kappa }\left (x  \right )=\lim_{x\to -\infty}g_{\kappa}(x)=0;
$$
$(iii)$  If $\kappa=1,$  then
$$
g_{\kappa}(x)\equiv e^{-e^{- \gamma}}>\frac{1}{2}.
$$
\end{pro}

\noindent {\bf Proof. }
Define a function
\begin{gather*}\label{b1}
h_{\kappa}(x):=(\kappa-1)x+\kappa \gamma,\ x\in \mathbb{R}.
\end{gather*}
Then
$$g_{\kappa}(x)=e^{-e^{-h_{\kappa}(x)}},\ x\in \mathbb{R}.$$
By taking the derivative of $g_{\kappa}(x)$, we have
\begin{gather*}\label{e1}
g^{'}_{\kappa}(x)=(\kappa-1)e^{-h_{\kappa}(x)}e^{-e^{-h_{\kappa}(x)}}.
\end{gather*}

If $\kappa <1,$ then $g^{'}_{\kappa}(x)< 0,$  which implies that $g_{\kappa}(x)$  is a strictly  decreasing function of $x.$  Then
$$
\inf_{x \in \mathbb{R} }g_{\kappa }\left (x  \right )=\lim_{x\to -\infty}g_{\kappa}(x)=0.
$$
If $\kappa> 1,$  then $g^{'}_{\kappa}(x)> 0.$  Thus, $g_{\kappa}(x)$  is a strictly increasing  function of $x$  and
$$
\inf_{x\in \mathbb{R}}g_{\kappa}(x)=\lim_{x \rightarrow -\infty}g_{\kappa}(x)=0.
$$
If $\kappa=1,$  then by (\ref{4.1-a}), we get that $g_{\kappa}(x)\equiv e^{-e^{- \gamma}}>\frac{1}{2}$. The proof is complete. \hfill\fbox

\subsection{Logistic distribution}\setcounter{equation}{0}

Let $X_{\mu,\beta}$  be a logistic random variable with parameters $\mu$ and $\beta$ $(\mu \in \mathbb{R}, \beta>0)$. By \cite[Chapter 21]{K}, we know that the distribution function of $X_{\mu,\beta}$ is given by
$$
P(X\leq x)=\frac{1}{1+e^{-\frac{x-\mu}{\beta}}}, \ x \in \mathbb{R},
$$
and the expectation of $X_{\mu,\beta}$ is $E[X]=\mu$. Then, for any given real number $\kappa>0,$  we have
$$
P(X_{\mu,\beta}\leq \kappa E[X_{\mu,\beta}])= \frac{1}{1+e^{-\frac{\kappa \mu-\mu}{\beta}}}.
$$

Let  $y:=\frac{\mu}{\beta}$.  Then
$$
P(X_{\mu,\beta}\leq \kappa E[X_{\mu,\beta}])= \frac{1}{1+e^{-(\kappa-1)y}}.
$$
Define a function
\begin{gather}\label{b2}
g_{\kappa}(y):=\frac{1}{1+e^{-(\kappa-1)y}},\  \  y\in \mathbb{R}.
\end{gather}

The main result of this section is
\begin{pro}\label{pro-5.2}
$(i)$ If $\kappa < 1$, then
\begin{eqnarray*}
\inf\limits_{y \in \mathbb{R} }g_{\kappa }\left (y \right )=\lim_{y\to \infty}g_{\kappa}(y)=0;
\end{eqnarray*}
$(ii)$ If $\kappa> 1$, then
$$
\inf_{y \in \mathbb{R} }g_{\kappa }\left ( y  \right )=\lim_{y\to -\infty}g_{\kappa}(y)=0;
$$
$(iii)$  If $\kappa=1,$  then
$$g_{\kappa}(y)\equiv \frac{1}{2}.$$

\end{pro}
{\bf Proof. } By the definition $g_{\kappa}(y)$  in \eqref{b2}, we have
$$g^{'}_{\kappa}(y)=\frac{(\kappa-1)e^{-(\kappa-1)y}}{(1+e^{-(\kappa-1)y})^2},\  \   y\in \mathbb{R}.$$
If $\kappa<1,$  then $g^{'}_{\kappa}(y)<0$, which implies that the function $g_{\kappa}(y)$  is strictly decreasing. Then
$$
\inf_{y\in \mathbb{R}} g_{\kappa}(y)=\lim_{y\rightarrow \infty}g_{\kappa}(y)=0.
$$
If $\kappa> 1,$ then $g^{'}_{\kappa}(y)> 0.$  Thus, $g_{\kappa}(y)$  is a strictly increasing function of $y$   and
$$
\inf_{y\in \mathbb{R}} g_{\kappa}(y)=\lim_{y\rightarrow -\infty}g_{\kappa}(y)=0.
$$
If $\kappa=1,$  then $g_{\kappa}(y)\equiv\frac{1}{2}$ by \eqref{b2}.  The proof is complete. \hfill\fbox

\bigskip

\noindent {\bf\large Acknowledgments}\quad  This work was supported by the National Natural Science Foundation of China (12171335, 12301603), the Science Development Project of Sichuan University (2020SCUNL201) and  the Scientific Foundation of Nanjing University of Posts and Telecommunications (NY221026).


\begin{thebibliography}{1234}

\bibitem{BPR23} L. Bababesi, L. Pratelli, P. Rigo, On the Chv\'{a}tal-Janson conjecture, Statis. Probab. Lett.  194 (2023) 109744.


\bibitem{Do18} B. Doerr,  An elementary analysis of the probability that a binomial random variable exceeds its expectation, Statis. Probab. Lett.  139 (2018) 67-74.



\bibitem{GM14} S. Greenberg, M. Mohri,  Tight lower bound on the probability of a binomial exceeding its expectation,   Statis. Probab. Lett. 86 (2014) 91-98.

\bibitem{GTH23} Z.-Y. Guo, Z.-Y. Tao, Z.-C. Hu, A study on the negative binomial distribution motivated by Chv\'{a}tal's theorem, arXiv: 2311.06877v1 (2023).


\bibitem{H63}  C. C. Heyde, On a property of the lognormal distribution,  J. R. Stat. Soc. Series B, 25 (1963) 392-393.

\bibitem{Ja21} S. Janson, On the probability that a binomial variable is at most its expectation,  Statis. Probab. Lett. 171 (2021) 109020.

\bibitem{K}K. Krishnamoorthy, Handbook of Statistical Distributions with Applications. CRC Press, New York, 2020.

\bibitem{LHZ23} C. Li, Z.-C. Hu,  Q.-Q, Zhou, A study on the Weibull and Pareto distributions motivated by Chv\'{a}tal's  theorem, Accepted by J. of Math. (PRC) (2023).

\bibitem{LXH23} F.-B. Li, K. Xu, Z.-C. Hu, A study on the Poisson, geometric and Pascal distributions motivated by Chv\'{a}tal's conjecture, Statis. Probab. Lett. 200 (2023) 109871.


\bibitem{S73}   F. W. Steutel,  Some recent results in infinite divisibility, Stoch. Proc. Appl. 1(1973) 125-143.

\bibitem{Su21} P. Sun,  Strictly unimodality of the probability that the binomial distribution is more than its expectation, Discret. Appl. Math.  301 (2021) 1-5.

\bibitem{SHS23} P. Sun, Z.-C. Hu, W. Sun, The extreme values of two probability functions for the Gamma distribution, arXiv: 2303.17487 (2023).


\end{thebibliography}
\end{document}